\documentclass{article}
\usepackage[utf8]{inputenc}
\usepackage{amsmath}
\usepackage{amsthm}
\usepackage{amssymb}
\usepackage{mathrsfs}

\title{Operator Lie Algebras of Rotations and Transformations in White Noise}

\author{
Wolfgang Bock\\
{\small Technomathematics Group}\\
{\small	University of Kaiserslautern}\\
{\small	P.\ O.\ Box 3049, 67653 Kaiserslautern, Germany}\\
{\small E-Mail:bock@mathemaik.uni-kl.de}\\[.3cm]
Janeth Canama\\
{\small MSU-IIT Iligan, Andres Bonifacio Avenue, Tibanga},\\ 
{\small	9200 Iligan City, Philippines}\\
{\small E-Mail:janeth.canama@g.msuiit.edu.ph}}

\date{\today}

\newtheorem{thm}{Theorem}[section]
\newtheorem{lem}[thm]{Lemma}

\newtheorem{rmk}[thm]{Remark}
\newtheorem{defn}[thm]{Definition}
\newtheorem{prop}[thm]{Proposition}
\newtheorem{cor}[thm]{Corollary}

\begin{document}

\maketitle

\begin{abstract}

The infinitesimal generator of a one-parameter subgroup of the infinite dimensional rotation group associated with the complex Gelfand triple
$
(E) \subset
L^2(E^*, \mu) \subset
(E)^*
$
is of the form
$$
R_\kappa
=
\int_{T\times T}
\kappa(s,t) (a_s^* a_t - a_t^* a_s) ds dt
$$
where $\kappa \in E \otimes E^*$
is a skew-symmetric distribution.
Hence
$R_\kappa$
 is twice the conservation operator
associated with a skew-symmetric operator $S$.
The Lie algebra containing $R_\kappa$, identity operator, annihilation operator, creation operator, number operator, 
(generalized) Gross Laplacian is discussed. We show that this Lie algebra is associated with the orbit of the skew-symmetric operator $S$.
\end{abstract}

\section{Introduction}
White noise calculus was initiated by Hida in 1975  \cite{Hida1975} and became a useful tool in many different areas, such as Mathematical Physics and Finance, see e.g.~ \cite{Hida1993,Obata1994,Kuo2,LetUsUse}. 

The root of white noise calculus is to switch a functional of Brownian motion 
$f(B(t); t \in \mathbb{R})$
with one of white noise
$\phi(\dot B(t); t\in \mathbb{R})$,
where $\dot B(t)$ is a time derivative of a Brownian motion $B(t)$.
We may thereby regard $\{{\dot B} (t)\}$
as a collection of infinitely many independent random variables and hence a coordinate system of an infinite dimensional space, \cite{Obata1994}.
Let $\mathcal{S}(\mathbb{R}^n)$ and
$\mathcal{S}'(\mathbb{R}^n)$ be the
Schwartz space consisting of rapidly decreasing
$C^\infty$-functions and the space of tempered distributions, respectively.
The mathematical framework of white noise calculus is based on an infinite dimensional analogue of Schwartz' distribution theory, where the roles of the Lebesgue measure on $\mathbb{R}^n$ and the Gelfand triple
$
\mathcal{S}(\mathbb{R}^n) \subset
L^2(\mathbb{R}^n) \subset
\mathcal{S}'(\mathbb{R}^n)
$
are played by the Gaussian measure $\mu$
on $E^*$ (the topological dual space of a nuclear space $E$)
and 
$(E) \subset 
(L^2) = L^2(E^*,\mu) 
\subset (E)^*$, respectively.
Furthermore, in white noise calculus the coordinate differential operators are given by 
$a_t$, where $t$ runs over a time parameter space $T$, 
\cite{Obata1992}.

The finite dimensional Laplacian $\Delta_n$ on $\mathbb{R}^n$ admits two expressions:
\begin{equation} \label{LaplaciansFinite}
\Delta_n
=
\sum_{i=1}^n \left( \frac{\partial}{\partial_{x_i}}\right)^2
=
- \sum_{i=1}^n
\left( \frac{\partial}{\partial_{x_i}}\right)^*
\frac{\partial}{\partial_{x_i}},
\end{equation}
when $\Delta_n$ acts on $\mathcal{S}(\mathbb{R}^n)$.
By virtue of a general theory established in
\cite{Hida1992},
the Gross Laplacian and the number operator are expressed as follows:
$$
\Delta_G = \int_T a_t^2 dt, \quad
N = \int_T a_t^* a_t dt.
$$
These are infinite dimensional analogues of the finite dimensional Laplacian.
However, unlike the finite dimensional case,
$\Delta_G$ and $N$ are completely different from each other \cite{Chung1998}.
It was Gross \cite{Gross1967} and Piech \cite{Piech1974} who
initiated the study of the Gross Laplacian and the number operator, as natural infinite-dimensional
analogues of a finite-dimensional Laplacian, in connection with the Cauchy problem in infinite-dimensional abstract Wiener space.
Based on white noise analysis, Kuo \cite{Kuo1986}
formulated the Gross Laplacian $\Delta_G$ and the number operator $N$ as continuous linear operators acting on white noise functionals. 
Hida, et.al. \cite{Hida1992} developed a general theory of operators acting on white noise functionals. 
As a particular case, they discussed infinite dimensional rotations in detail.
It was proved that $\Delta_G$ and $N$ are in essence the only operators which are rotation-invariant \cite{Obata1992}.
Moreover, characterization theorems for 
$\Delta_G$, $N$  and the Euler operator $\Delta_G + N$ were given in \cite{He1996}.
Chung, et. al. \cite{Chung1998} studied their generalizations 
$N(f)$ and $\Delta_G(g)$ which are called second order differential operators of diagonal type.
A generalized Gross Laplacian $\Delta_G(K)$, called the $K$-Gross Laplacian, was introduced in \cite{Chung1997}. 
Later on, the $K$-Gross Laplacian, the second quantization
and the differential second quantization were studied within the framework of nuclear algebras of entire functions 
\cite{Barhoumi2010}.

In white noise analysis, Lie algebraic approach in studying the infinite dimensional Laplacians $N$ and $\Delta_G$ have been  discussed in \cite{Chung2000}.
In 1998, Chung and Ji \cite{Chung1998T}
constructed a two-parameter transformation group
which is the two-dimensional Lie group associated with the Lie algebra $\mathbb{C} \Delta_G + \mathbb{C}N$. Obata \cite{Obata1995L} discussed all possible two-dimensional complex Lie algebras containing $N$ and constructed the associated Lie groups.
Chung and Chung \cite{Chung1996} studied the Lie algebras of Wick derivations on $(E)^*$.
Furthermore, Hida, et.al. \cite{Hida1993} obtained
a five-dimensional complex Lie algebra
$\mathfrak{h}$
 generated by the identity operator $Id$, $\Delta_G$, $N$, and
the infinitesimal generators of some differentiation and multiplication operators. 
Later on, Chung and Ji \cite{Chung2000}
explicitly constructed a Lie group associated with $\mathfrak{h}$.
They showed that this Lie algebra is spanned by $Id, a(\zeta), a^*(\zeta), N$ and 
$\Delta_G$, where 
$a(\zeta)$ and  $a^*(\zeta)$ are the annihilation and creation operators, respectively.
In 2016, Ji and Sinha \cite{Ji2016},
studied a class of fundamental quantum stochastic processes induced by the generators of a six dimensional non-solvable 
Lie $*$-algebra consisting of all linear combinations of the
identity operator, generalized Gross Laplacian and its adjoint, annihilation operator, creation operator, and conservation operator.

The paper is organized as follows:
In Section 2 we assemble standard notations used in white noise calculus and in Section 3 we discuss the infinite dimensional rotation group.
In Section 4 we review the most basic notions in quantum white noise calculus.
In Section 5 we survey some known commutation relations of white noise operators.
Section 6 is devoted to a study of the Lie algebra generated by the identity operator,  annihilation operator, 
creation operator, Gross Laplacian, generalized Gross Laplacian, number operator and rotation (conservation) operator.
We show that this Lie algebra is associated with the orbits of a skew-symmetric operator.
\section{White noise calculus}
If $\mathfrak{X}$ is a locally convex space over $\mathbb{R}$, its complexification is denoted by 
$\mathfrak{X}_\mathbb{C}$ \cite{Obata1994}.
The following construction of the complex Gelfand triple
is lifted from \cite{Chung1998}.
The whole discussion is based on the special choice of a real Gelfand triple
\begin{equation} \label{RealGelfand}
E = \mathcal{S}(\mathbb{R}) \subset
H = L^2(\mathbb{R}, dt) \subset
E^* =  \mathcal{S}'(\mathbb{R}).
\end{equation}
However, $\mathbb{R}$ can be replaced with $\mathbb{R}^n$ with no essential change (sometimes more interesting for applications).
It is noteworthy that (\ref{RealGelfand})
is constructed from the differential operator
$A = 1+t^2 -\frac{d^2}{dt^2}$.
In fact, $E =  \mathcal{S}(\mathbb{R})$
is identified (up to null functions) with the space of functions $\xi \in H$ such that
$|\xi|_p = |A^p \xi|_0 < \infty$ 
for any $p \in \mathbb{R}$, where
$|\cdot|_0$ stands for the norm of $H$, and the topology of $E$ is given by the norms
$|\cdot |_p, p \in \mathbb{R}$.
Since $A$ is a positive self-adjoint operator with 
Hilbert-Schmidt inverse,
$E$ becomes a countable Hilbert nuclear space.
By definition $E^*$ is the strong dual space of $E$.
The canonical bilinear form on $E^* \times E$
and the real inner product of $H$ are denoted by the same symbol $\langle \cdot, \cdot \rangle$ because they are consistent.
The  Gaussian measure $\mu$ is by definition a unique probability measure on $E^*$ of which characteristic function is
$$
\exp \left( -\frac{1}{2} |\xi|^2_0 \right)
= \int_{E^*} e^{i\langle x, \xi\rangle} \mu(dx), \quad
\xi \in E.
$$
The probability space $(E^*, \mu)$ is called the
 white noise space or the  Gaussian space.
With each $\xi \in E_\mathbb{C}$, we associate a function on $E^*$ defined by
$$
\phi_\xi(x) = 
\exp
\left(
\langle x, \xi \rangle - \frac{1}{2} 
\langle \xi, \xi \rangle
\right ), 
\quad x\in E^*,
$$
which is called an  exponential vector.
The correspondence
$$
\phi_\xi \longleftrightarrow
\left(
1, \frac{\xi}{1!}, \frac{\xi^{\otimes2}}{2!},
\frac{\xi^{\otimes 3}}{3!}, \dots
\right), 
\quad \xi \in E_{\mathbb{C}},
$$
is uniquely extended to a unitary isomorphism between
$L^2(E^*, \mu)$ and the Boson Fock space over
$H_\mathbb{C}$, denoted by $\Gamma(H_\mathbb{C})$,
which is the celebrated
 Wiener-Ito-Segal isomorphism.
If $\phi \in L^2(E^*, \mu)$ and 
$(f_n) \in \Gamma(H_\mathbb{C})$ are related,
we write $\phi \sim (f_n)$ simply. 
In that case,
$$
||\Phi||^2_0  =
\int_{E^*} |\phi(x)|^2 \mu(dx) =
\int_{n=0}^\infty n! |f_n|^2_0
.
$$
For any $p \in \mathbb{R}$ we put
$$
|| \phi ||^2_p = 
\int_{n=0}^\infty n! |f_n|^2_p = 
\int_{n=0}^\infty n! |(A^{\otimes n})^p|^2_0, \quad
\phi \sim (f_n)
.
$$
Let $(E)$ be the subspace of functions 
$\phi \in L^2(E^*, \mu)$ such that
$\| \phi \|_p < \infty$ for all $p$.
Then $(E)$ becomes a nuclear Frechet space
with the defining seminorms $\| \cdot \|_p, p \in \mathbb{R}$.
The dual space $(E)^*$ consists of all elements
$\Phi \sim (F_n)$ such that
$F_n \in (E^{\otimes n}_\mathbb{C})^*_{sym}$
and $|| \Phi ||_{-p} < \infty$ for some $p \ge 0$.
We thereby obtain a complex Gelfand triple:
$$
(E) \subset L^2(E^*, \mu) \subset (E)^*.
$$
Elements in $(E)$ and $(E)^*$ are called
a  test (white noise) function and
a  generalized (white noise) function, respectively.

A continuous linear operator from $(E)$ into $(E)^*$ is called a
white noise operator. The space of white noise operators is denoted by $\mathcal{L}((E),(E)^*)$ and is equipped with the bounded convergence topology. It is noted that
$\mathcal{L}((E),(E))$ is a subspace of 
$\mathcal{L}((E),(E)^*)$.
With each $y \in E^*_\mathbb{C}$ we may associate an 
 annihilation operator $D_y \in L((E), (E))$ 
which is uniquely determined by 
$$
D_y \phi_\xi = \langle y, \xi \rangle \phi_\xi,
 \
\xi \in E_\mathbb{C}.
$$
Since $\delta_t \in E^*$ for any $t \in \mathbb{R}$,
$$
a_t = D_{\delta_t}, \quad t\in \mathbb{R},
$$
belongs to $L((E),(E))$.
This is called the  annihilation operator
at a point $t \in \mathbb{R}$.
The  creation operator at a point
is by definition the adjoint 
$a_t^* \in L((E)^*, (E)^*)$.
It is known that $a_t$ is a differential operator along the direction $\delta_t$, namely
$$
a_t \phi(x) =
\lim_{\theta \to 0}
\frac{\phi(x+\theta \delta_t)-\phi(x)}{\theta}, \quad
\phi \in (E), t \in \mathbb{R}, x \in E^*.
$$

\begin{defn} \cite{Obata1995D}  
With each 
$\kappa \in (E_{\mathbb{C}}^{\otimes(l+m)})^*$
we may associate an  integral kernel operator
whose formal expression is given by
$$
\Xi_{l,m}(\kappa)
=
\int_{T^{l+m}}
\kappa(s_1, \cdots, s_l, t_1, \cdots, t_m)
a_{s_1}^* \cdots a_{s_l}^*
a_{t_1} \cdots a_{t_m}
ds_1 \cdots ds_l
dt_1 \cdots dt_m,
$$
where $\kappa$ is called the 
 kernel distribution.
\end{defn}

\begin{prop} \cite{Obata1994}
\label{integral-kernel-char}
Let $\phi \in (E)$ be given with Wiener-Ito expansion:
$$
\phi(x) = \sum_{n=0}^\infty 
\langle 
: x^{\otimes n}:, f_n 
\rangle.
$$
Then, for $\kappa \in (E_\mathbb{C}^{\otimes (l+m)})^*$
we have
$$
\Xi_{l,m}(\kappa) \phi(x) = 
\sum_{n=0}^\infty \frac{(n+m)!}{n!}
\langle
:x^{\otimes (l+n)}:, \kappa \otimes_m f_{n+m}
\rangle.
$$
\end{prop}

\begin{lem} \cite{Chung1998} \label{Char-Integral-Kernel}
Let $\kappa \in (E_\mathbb{C}^{\otimes (l+m)})^*$.
Then $\Xi_{l,m}(\kappa) \in \mathcal{L}((E),(E))$ if and only if
$\kappa \in (E_\mathbb{C}^{\otimes l}) \otimes (E_\mathbb{C}^{\otimes m})^*$. In particular, 
$\Xi_{0,m} \in \mathcal{L}((E),(E))$ for any
$\kappa \in (E_\mathbb{C}^{\otimes m})^*$.
\end{lem}


\begin{prop} \cite{Obata1994}
For any $T \in \mathcal{L}(E_\mathbb{C}, E_\mathbb{C})$,
there exists a unique operator $\Gamma(T) \in
\mathcal{L}((E),(E))$ such that
$
\Gamma(T) \phi_\xi = \phi_{T\xi}, \ 
\xi \in E_\mathbb{C}.
$ 
Moreover, for $\phi \in (E)$ given with
Wiener-Ito expansion:
$$
\phi(x) = \sum_{n=0}^\infty 
\langle : x^{\otimes n}:, f_n \rangle
$$
it holds that
$$
\Gamma(T) \phi(x) = 
\sum_{n=0}^\infty
\langle : x^{\otimes n}:, T^{\otimes n} f_n \rangle
$$
\end{prop}
In general, for $T \in \mathcal{L}(E_\mathbb{C}, E_\mathbb{C})$ we define an operator $d\Gamma(T)$ on $(E)$.
Suppose $\phi \in (E)$ is given as
$$
\phi(x) = \sum_{n=0}^\infty
\langle
:x^{\otimes n}:, f_n
\rangle,
\quad 
x\in E^*,
$$
as usual. Then we put
$$
d\Gamma(T) \phi(x) =
\sum_{n=0}^\infty 
\langle
:x^{\otimes n}:, 
\gamma_n(T) f_n
\rangle,
$$
where
$$
\begin{cases}
\gamma_n(T) & =
\sum_{k=0}^{n-1} I^{\otimes k} \otimes T \otimes
I^{\otimes (n-1-k)}, \ n \ge 1,\\
\gamma_0(T) &= 0.
\end{cases}
$$

\begin{prop} \cite{Obata1994}
$d\Gamma(T) \in \mathcal{L}((E), (E))$
for any $T \in \mathcal{L}(E_\mathbb{C}, E_\mathbb{C})$.
\end{prop}

\begin{defn} \cite{Obata1994}
$\Gamma(T)$ and $d\Gamma(T)$ are called
the second quantization and differential second quantization of $T$, respectively.
\end{defn}

If $K$ and $\kappa \in E_\mathbb{C} \otimes E_\mathbb{C}^*$
are related by the kernel theorem, i.e.
\begin{equation} \label{kernel-theorem}
\langle K\xi, \eta \rangle =
\langle \kappa, \eta \otimes \xi \rangle, \quad
\xi, \eta \in E_\mathbb{C},
\end{equation}
then $d\Gamma(K) = \Xi_{1,1}(\kappa)$.

\subsection{Second order differential operators of diagonal type}
Second order differential operators of diagonal type and their generalizations were discussed in \cite{Chung1998}.
Let $\tau \in (E_\mathbb{C} \otimes E_\mathbb{C})^*$
be defined by
$$
\langle \tau, \eta \otimes \xi \rangle =
\langle \xi, \eta \rangle, \quad
\xi, \eta \in E_\mathbb{C}.
$$
In fact, $\tau \in E_\mathbb{C} \otimes E_\mathbb{C}^*$
since $\tau$ corresponds to the identity operator under the canonical isomorphism 
$E_\mathbb{C} \otimes E_\mathbb{C}^* \cong
\mathcal{L}(E_\mathbb{C}, E_\mathbb{C})$.
Thus, by Lemma \ref{Char-Integral-Kernel},
$$
N = \Xi_{1,1}(\tau) = 
\int_{\mathbb{R}^2} \tau(s,t) a_s^* a_t ds dt
=
\left(
\int_\mathbb{R}
a_t^* a_t dt \text{ for simplicity}
\right)
$$
belongs to $\mathcal{L}((E),(E))$.
This is called the number operator.
On the other hand,
$$
\Delta_G = \Xi_{0,2}(\tau)
=
\int_{\mathbb{R}^2}
\tau(t_1, t_2) a_{t_1} a_{t_2} dt_1 dt_2
=
\left(
\int_{\mathbb{R}} a_t^2 dt \text{ for simplicity}
\right)
$$
also belongs to $\mathcal{L}((E),(E))$ 
and is called the Gross Laplacian. 

\subsection{Convolution}
Suppose $S_1 \in \mathcal{L}(E_\mathbb{C}, E_\mathbb{C})$ and
$S_2 \in \mathcal{L}(E_\mathbb{C}, E_\mathbb{C}^*)$.
Let $f_1 \in E_\mathbb{C} \otimes E_\mathbb{C}^*$ and
$f_2 \in (E_\mathbb{C} \otimes E_\mathbb{C})^*$ be the corresponding elements, 
respectively, see (\ref{kernel-theorem}).
Then we denote by $f_2 *   f_1$ the element of
$(E_\mathbb{C} \otimes E_\mathbb{C})^*$ corresponding to 
$S_2 S_1 \in \mathcal{L}(E_\mathbb{C}, E_\mathbb{C}^*)$.
It is noted in \cite{Chung1998} that
$$
f_2 *   f_1 (s,t) = 
\int_{\mathbb{R}} 
f_2(s,u) f_1(u,t) du
$$
in a generalized sense.

\begin{lem} \cite{Chung1998}
\label{convolution}
The convolution $(f_2, f_1) \mapsto f_2 *   f_1$
gives a separately continuous bilinear map from
$(E_\mathbb{C} \otimes E_\mathbb{C})^* \times
(E_\mathbb{C} \otimes E_\mathbb{C}^*)$ into
$(E_\mathbb{C} \otimes E_\mathbb{C})^*$, and
from 
$(E_\mathbb{C} \otimes E_\mathbb{C}^*) \times
(E_\mathbb{C} \otimes E_\mathbb{C}^*)$ into
$E_\mathbb{C} \otimes E_\mathbb{C}^*$.
\end{lem}

\section{Infinite dimensional rotation group}
Our discussion is based on a Gelfand triple
$E \subset H \subset E^*$.
We denote by $GL(E)$ the group of all linear homeomorphisms from $E$ onto itself.
Then $GL(E) \subset L(E,E)$.
We now set
$$
O(E;H) =
\{
g \in GL(E) :
|g\xi|_0 = |\xi|_0 \text{ for all }
\xi \in E
\},
$$
which is a subgroup of $GL(E)$.
The group $O(E;H)$ is called the
 infinite dimensional rotation group
(associated with the Gelfand triple 
$E \subset H \subset E^*$).
The Gaussian measure $\mu$
is invariant under $O(E;H)$ \cite{Obata1994}.
Similarly we put
$$
U((E);(L^2)) =
\{
U \in GL((E)) :
\| U\phi\|_0 = 
\| \phi \|_0,
\text{ for all }
\phi \in (E)
\}.
$$
With each $g \in O(E;H)$ we associate its second
quantization $\Gamma(g)$.

\begin{defn} \cite{Obata1994}
We say that a continuous operator from $(E)$ into $(E)^*$ 
is rotation-invariant if
\begin{equation} \label{rotation-invariance}
\Gamma(g) \Xi \Gamma(g) = \Xi \quad
\text{ for all } \quad
g \in O(E;H).
\end{equation}
\end{defn}

\begin{rmk} \cite{Obata1994} 
The condition (\ref{rotation-invariance}) for
$\Xi \in \mathcal{L}((E), (E))$ is equivalent to the following
$
[\Gamma(g), \Xi] = 0 \quad
\text{ for all} \quad
g \in O(E;H).
$
\end{rmk}
In the following discussion let $\mathfrak{X}$ be a nuclear Frechet space with defining Hilbertian seminorms 
$
\{ \| \cdot \| \}_{\alpha \in A}
$
taking $\mathfrak{X} = E$ or
$\mathfrak{X} = (E)$ into consideration.

\begin{defn} \cite{Obata1994}
A one-parameter subgroup 
$\{g_\theta\}_{\theta \in \mathbb{R}}
 \subset GL(\mathfrak{X})$
is called differentiable if
$$
\lim_{\theta \to 0} \frac{g_\theta \xi - \xi}{\theta}
$$
converges in $\mathfrak{X}$ for any $\xi \in \mathfrak{X}$.
In that case a linear operator $X$ from $\mathfrak{X}$ into itself is defined by
$$
X\xi = 
\lim_{\theta \to 0} \frac{g_\theta \xi - \xi}{\theta},
\quad \xi \in \mathfrak{X}.
$$
This operator $X$ is called the infinitesimal generator of the diffentiable one-parameter subgroup
$\{g_\theta\}_{\theta \in \mathbb{R}} \subset
GL(\mathfrak{X})$.
\end{defn}

\begin{prop} \cite{Obata1994}
Let $\{g_\theta\}_{\theta \in \mathbb{R}} \subset
GL(\mathfrak{X})$ be a differentiable one-parameter
subgroup. Then its infinitesimal generator $X$ is always continuous, i.e., $X \in \mathcal{L}(\mathfrak{X}, \mathfrak{X})$. 
\end{prop}

\begin{defn} \cite{Obata1994}
A differentiable one-parameter subgroup
$\{g_\theta \}_{\theta \in \mathbb{R}} \subset
GL(\mathfrak{X})$ with infinitesimal generator $X$
is called regular if for any $\alpha \in A$
there exists $\beta \in A$ such that
$$
\lim_{\theta \to 0}
\sup_{\| \xi \|_\beta \le 1}
\left
\|
\frac{g \theta \xi - \xi}{\theta} - X\xi
\right
\|_\alpha
= 
0.
$$
\end{defn}

\begin{thm} \cite{Obata1994}
For $y \in E^*$, let $D_y$ and $T_y$ be the differential and translation operators, respectively.
Then $\{T_{\theta y}\}_{\theta \in \mathbb{R}}$ is a regular
one-parameter subgroup of $GL((E))$ with infinitesimal generator $D_y$.
\end{thm}

\begin{rmk} \cite{Obata1994} \label{skew-symmetric}
The infinitesimal generator $X$ of a regular one-parameter subgroup
$\{ g_\theta\}_{\theta \in \mathbb{R}}$
is skew-symmetric in the sense that
\begin{equation} \label{skew-symmetric-gen}
\langle X \xi, \eta \rangle = 
-\langle \xi, X \eta \rangle, 
\quad \xi, \eta \in E.
\end{equation}
Let $X \in \mathcal{L}(E_\mathbb{C}, E_\mathbb{C})$
be skew-symmetric in the sense of 
(\ref{skew-symmetric-gen}).
By the canonical isomorphism
$E_\mathbb{C} \otimes E_\mathbb{C}^*
\cong \mathcal{L}(E_\mathbb{C}, E_\mathbb{C})$
there exists $\kappa \in E_\mathbb{C} \otimes E_\mathbb{C}^*$
such that
$$
\langle \kappa, \eta \otimes \xi \rangle
= \frac{1}{2} \langle \eta, X \xi \rangle,
\quad \xi, \eta \in E_\mathbb{C}.
$$
Since $X$ is skew-symmetric, we have
\begin{equation} \label{skew-symmetric-dist}
\langle \kappa, \eta \otimes \xi \rangle =
\frac{1}{2}  \langle \eta, X \xi \rangle =
-\frac{1}{2}  \langle \xi, X\eta \rangle =
- \langle \kappa, \xi \otimes \eta \rangle.
\end{equation}
\end{rmk}

\begin{thm}{\cite{Obata1994}} \label{rotation-operator}
Let $\{g_\theta\}_{\theta \in \mathbb{R}}$
be a regular one-parameter subgroup of
$O(E;H)$ with infinitesimal generator $X$.
Then 
 $\{ \Gamma (g_\theta) \}_{\theta \in \mathbb{R}}$
is a regular one-parameter subgroup of
$U((E);(L^2))$ with infinitesimal generator
$d\Gamma(X)$. Moreover, there exists a skew-symmetric
distribution $\kappa \in E \otimes E^*$ such that
\begin{equation} \label{rotation-integral-kernel}
d\Gamma(X) = \int_{T\times T}
\kappa(s,t) (a_s^* a_t - a_t^* a_s) ds dt
= 2\Xi_{1,1}(\kappa).
\end{equation}
\end{thm}

\section{Quantum white noise calculus}


For $f \in E^*$ we define white noise operators:
\begin{equation}
a(f) = \Xi_{0,1}(f) = 
\int_T f(t) a_t dt, 
\qquad
a^*(f) = \Xi_{1,0}(f) =
\int_T f(t) a_t^* dt,
\end{equation}
which are called respectively the annihilation and creation operators associated with $f$.
If $\zeta \in E$, then $a(\zeta)$ extends to a continuous linear operator from $(E)^*$ into itself 
(denoted by the same symbol) and
$a^*(\zeta)$ (restricted to $(E)$) is a continuous linear operator from $(E)$ into itself \cite{Ji2007}.

\begin{lem} \cite{Ji2010Q}
For $\zeta \in E$, both $a(\zeta)$ and $a^*(\zeta)$
belong to 
$\mathcal{L}((E),(E)) \cap 
\mathcal{L}((E)^*, (E)^*)$.
\end{lem}

Thus, for any white noise operator 
$\Xi \in \mathcal{L}((E), (E)^*)$ and $\zeta \in E$, the commutators
$$
[a(\zeta), \Xi ] = a(\zeta) \Xi - \Xi a(\zeta), \quad
-[a^*(\zeta), \Xi] = \Xi a^*(\zeta) - a^*(\zeta) \Xi,
$$
are well-defined white noise operators, i.e.,
belongs to $\Xi \in \mathcal{L}((E), (E)^*)$.
We define
$$
D^+_\zeta \Xi = [a(\zeta), \Xi], \quad
D^-_\zeta \Xi = - [a^*(\zeta), \Xi].
$$
\begin{defn} \cite{Ji2007}
$D^+_\zeta \Xi$ and $D^-_\zeta \Xi$ are respectively called the creation derivative and annihilation derivative of $\Xi$, 
and both together the quantum white noise derivatives
(qwn-derivatives for brevity) of $\Xi$.
\end{defn}

\begin{thm} \cite{Ji2007}
The map
$$
E \times \mathcal{L}((E),(E)^*) \to \mathcal{L}((E),(E)^*), \quad
(\zeta, \Xi) \mapsto D^{\pm}_\zeta \Xi
$$
is continous bilinear.
\end{thm}

\begin{cor} \cite{Ji2007}
For each $\zeta  \in E$, 
the qwn-differential operator $D^{\pm}_\zeta$
is a continuous operator from
$\mathcal{L}((E),(E)^*)$ into itself.
\end{cor}
The quantum white noise derivatives of the generalized Gross Laplacian and conservation operator are given in 
\cite{Ji2010, Ji2010Q}.
For each $S \in \mathcal{L}(E, E^*)$, by the kernel theorem there exists a unique $\tau_S \in (E \otimes E)^*$ such that
$$
\langle \tau_S, \eta \otimes \xi \rangle =
\langle S \xi, \eta \rangle , 
\quad \xi, \eta \in E.
$$
The integral kernel operator
$$
\Delta_G(S) = \Xi_{0,2}(\tau_S) =
\int_{T\times T} \tau_S(s,t) a_s a_t ds dt
$$
is called the generalized Gross Laplacian associated with $S$.
Note that 
$\Delta_G (S)  \in \mathcal{L}((E),(E))$.
The usual Gross Laplacian is 
$\Delta_G = \Delta_G(I)$.
The integral kernel operator 
$$
\Lambda(S) = \Xi_{1,1}(\tau_S) = 
\int_{T \times T} \tau_S (s,t) a_s^* a_t ds dt
$$
is called the conservation operator associated with $S$.
In general, $\Lambda(S) \in \mathcal{L}((E),(E)^*)$.
Note that $N = \Lambda(I)$ is the number operator.

\begin{lem} \cite{Ji2010}
For $S \in \mathcal{L}(E, E^*)$ and $\zeta \in E$, we have
\begin{equation} \label{Qwn-derivatives-Gross}
D_\zeta^+ \Delta_G(S) = 0, \qquad 
D_\zeta^- \Delta_G(S) = a(S\zeta) + a(S^* \zeta),
\end{equation}
\begin{equation} \label{Qwn-derivatives-Number}
D_\zeta^+ \Lambda(S) = a(S^*\zeta), \qquad \qquad
D_\zeta^- \Lambda(S) = a^*(S\zeta).
\end{equation}
\end{lem}
\section{Commutation relations}
We have the so-called 
canonical commutation relation:
\begin{equation} \label{CCR}
[a_s, a_t] = 0, \quad
[a_s^*, a_t^*] = 0, \quad
[a_s, a_t^*] = \delta_s(t) Id, \quad
s,t \in \mathbb{R},
\end{equation}
where the last relation is understood in a generalized sense.

\begin{prop} \cite{Obata1994} \label{generalizedCCR}
For $\xi \in E_\mathbb{C}$ and
$y \in E^*_\mathbb{C}$, it holds that
$$
[
\Xi_{0,1}(y), 
\Xi_{1,0}(\xi)
]
=
\langle y, \xi \rangle Id
$$
where $Id$ is the identity operator on $(E)$.
\end{prop}

\begin{thm} \cite{Chung2000}  \label{ManyRelations}
For $\zeta \in E_\mathbb{C}$, we have the following commutation relations:
\begin{enumerate}
	\item $[a(\zeta), a^*(\zeta)] =
	\langle \zeta, \zeta \rangle Id$,
	\item $[a(\zeta), N] = a(\zeta)$,
	\item $[a(\zeta), \Delta_G] = 0$,
	\item $[a^*(\zeta), N] = -a^*(\zeta)$,
	\item $[a^*(\zeta), \Delta_G] = -2a(\zeta)$,
	\item $[\Delta_G, N] = 2\Delta_G$.
	\item $
[N, \Xi_{0,m} (\kappa)] = 
-m \Xi_{0,m} (\kappa), \quad
\kappa \in (E^{\otimes m}_\mathbb{C})^*.
$
\end{enumerate}
\end{thm}

\begin{lem} \cite{Chung1998}
\label{CCR-Number-operator}
It holds that
$$
[\Xi_{1,1}(f_1), \Xi_{1,1}(f_2)] = 
\Xi_{1,1}(f_1 * f_2 - f_2  * f_1), \qquad
f_1, f_2 \in E_\mathbb{C} \otimes E_\mathbb{C}^*.
$$
\end{lem}

\begin{rmk} \cite{Chung1998} \label{annihilation-commute}
In general, $\Xi_{0,2}(\kappa)$ with 
$\kappa \in (E_\mathbb{C} \otimes E_\mathbb{C})^*$ involves only annihilation operators, and so they commute each other.
\end{rmk}

\begin{thm} \cite{Chung2000}
\label{Lie-algebra-base}
For each nonzero $\zeta \in E_\mathbb{C}$, let
$\mathfrak{h} = 
\langle Id, a(\zeta), a^*(\zeta), N, \Delta_G \rangle.
$
Then $\mathfrak{h}$ is a five-dimensional non-nilpotent solvable complex Lie algebra.
\end{thm}

\begin{thm} \cite{Chung2000}
Let $\mathfrak{g} = 
\langle
N, \Xi_{0,m_1} (\kappa_1), \dots, \Xi_{0,m_n} (\kappa_n)
\rangle$
be the complex vector space spanned by $N$ and
$\Xi_{0,m_i} (\kappa_i),  \ i =1, \dots n$,
where for each \ $i = 1, \dots, n, 
\kappa \in (E_\mathbb{C}^{\otimes m_i})^*$.
Then $\mathfrak{g}$ is a $(n+1)$-dimensional non-nilpotent solvable complex Lie algebra.
\end{thm}

\begin{thm} \cite{Ji2016}
\label{Fixed-point}
Let $K \in \mathcal{L}(H, H)$ be a symmetric
Hilbert-Schmidt operator,
$L \in \mathcal{L}(H, H)$ a self-adjoint operator, 
and $\zeta \in H$ satisfying that
$$
KL = \overline{L} K = K, \quad
\overline{K} K = L, \quad
\overline{K\zeta} = \zeta, \quad
L\zeta = \zeta.
$$
Then, the linear span of
$Id, a(\zeta), a^*(\zeta), \Lambda(L), 
\Delta_G(K), \Delta^*_G(K)$
becomes a six dimensional non-solvable Lie 
$*$-algebra.
\end{thm}

\begin{rmk} \cite{Ji2016}
In Theorem \ref{Fixed-point},
the operators $K, L \in \mathcal{L}(H,H)$ 
satisfy that 
$\overline{L} K = K$ and 
$\overline{K} K = L$.
Therefore, we see that
$L \overline{K} K = \overline{K} K$ and so
$L^2 = L$.
Since $L$ is a self-adjoint operator, 
$L$ is a projection.
\end{rmk}

\section{Main Results}
Motivated by the result in Theorem \ref{rotation-operator},
given a skew-symmetric distribution $\kappa \in E \otimes E^*$,
we will study an operator whose expression is given by
$$
R_\kappa =
\int_{T\times T}
\kappa(s,t) (a_s^* a_t - a_t^* a_s) ds dt
=
2\Xi_{1,1}(\kappa)
.
$$
We enlarge the Lie algebra described in Theorem \ref{Lie-algebra-base}
by adding the operator
$\Xi_{1,1}(\kappa)$ and investigating their commutators.
\subsection{Notations}
We need notations. For any white noise operator
$\Xi \in \mathcal{L}((E), (E)^*)$ and $\zeta \in E$,
set
$$
D^{0+}_\zeta \Xi = D^{+}_\zeta \Xi, \qquad
D^{0-}_\zeta \Xi = D^{-}_\zeta \Xi . 
$$
Moreover, for $k =1, 2, \dots$, define
$$
D^{k+}_\zeta \Xi = [ D^{(k-1)+}_\zeta \Xi, \Xi], \qquad
D^{k-}_\zeta \Xi = -[D^{(k-1)-}_\zeta \Xi, \Xi]. 
$$

\subsection{Quantum white noise derivatives of an integral kernel operator}
\begin{lem}
\label{Qwn-derivatives-Xi}
Let $\zeta \in E$ and $S \in \mathcal{L}(E,E)$ be skew-symmetric
in the sense of (\ref{skew-symmetric-gen}).
Then
for $k = 0, 1, 2 \dots$, we have
\begin{align*}
D^{k-}_\zeta \Lambda(S) &= a^*(S^{k+1} \zeta), \qquad
D^{k+}_\zeta \Lambda(S)  = (-1)^{k+1}  a(S^{k+1} \zeta), \\
D^{k-}_\zeta \Delta_G(S) &= 0, 
\qquad \qquad \qquad
D^{k+}_\zeta \Delta_G (S) = 0. 
\end{align*}
\end{lem}
\noindent
\emph{Proof:}
Let $\Xi = \Lambda(S) = \Xi_{1,1}(\tau_S)$.
Since $S$ is skew-symmetric, 
then $\tau_S$ is a skew-symmetric distribution. Hence,
\begin{align*}
a(S^*\zeta) &=
\int_T S^*\zeta(t)  a_t dt \\
&=
\int_{T\times T} \tau_S(u,t) \zeta(u) a_t dt du\\
&= 
-\int_{T\times T} \tau_S(t,u) \zeta(u) a_t dt du\\
&=
-a(S\zeta).
\end{align*}
Applying (\ref{Qwn-derivatives-Number}), we have
\begin{equation}
\label{Dminus}
D^-_\zeta \Xi 
= -[a^*(\zeta), \Xi ]
=  a^*(S\zeta),
\end{equation}
\begin{equation}
\label{Dplus}
D^+_\zeta \Xi 
= [a(\zeta), \Xi ]
= a(S^* \zeta) = - a(S\zeta).
\end{equation}
Since $S \in \mathcal{L} (E,E)$, it follows that $S\zeta \in E$. Then applying (\ref{Dminus}) we have
\begin{align*}
D^{1-}_\zeta \Xi 
& = - [ D^{-}_\zeta \Xi, \Xi]
= -[ a^* (S \zeta)  ,\Xi]
= a^*(S^2 \zeta),\\
D^{2-}_\zeta \Xi 
& = - [ D^{1-}_\zeta \Xi, \Xi]
= -[ a^*(S^2 \zeta) , \Xi]
= a^*(S^3 \zeta),\\
& \vdots\\
D^{k-}_\zeta \Xi 
& = - [ D^{k-}_\zeta \Xi, \Xi]
= -[ a^*(S^k \zeta) , \Xi]
= a^*(S^{k+1} \zeta).
\end{align*}
Moreover, applying (\ref{Dplus}) we have
\begin{align*}
D^{1+}_\zeta \Xi 
& =  [ D^{+}_\zeta \Xi, \Xi]
= [-a ( S\zeta)  ,\Xi]
= a(S^2 \zeta),\\
D^{2+}_\zeta \Xi 
& =  [ D^{1+}_\zeta \Xi, \Xi]
= [ a(S^2 \zeta), \Xi]
= -a(S^3 \zeta),\\
& \vdots
\\
D^{k+}_\zeta \Xi & 
=  [ D^{(k-1)+}_\zeta \Xi, \Xi]
= [ (-1)^{k} a(S^k \zeta ) , \Xi]
= (-1)^{k+1} a(S^{k+1} \zeta). 
\end{align*}
Furthermore,
applying  (\ref{Qwn-derivatives-Gross}), we have
$$ 
D^-_\zeta \Delta_G(S)
= a(S\zeta) + a(S^*\zeta) 
= 0,
\qquad
D^+_\zeta \Delta_G(S)
= 0. \ \Box
$$

\begin{rmk}
From the known commutator identity 
(see, e.g. \cite{Lavrov2014})
$$
[A, BC] = [A,B]C + B[A,C]
$$
we can derive the following identity
\begin{equation} \label{comm-identity}
[AB,CD] = 
A[B,C]D + [A,C]BD + CA[B,D] + C[A,D]B.
\end{equation}
\end{rmk}

\begin{lem}
\label{Commutations3}
Let $\lambda, \kappa \in E_\mathbb{C} \otimes  E_\mathbb{C}^*$. 
 Then 
$$[\Xi_{0,2}(\lambda), \Xi_{1,1}(\kappa)] = 
\Xi_{0,2}(\lambda *   \kappa) 
+ \Xi_{0,2}(\lambda^* *   \kappa),
$$
where $\lambda^*(s,t) = \lambda(t,s)$ for all $s,t \in T$.
	\end{lem}
\noindent
\emph{Proof:}
Applying (\ref{comm-identity})
and (\ref{CCR}), we have
\begin{align*}
[a_s a_t, a_u^* a_v] 
& =
a_s [a_t, a_u^*] a_v + 
[a_s, a_u^*] a_t a_v +
a_u^*a_s [a_t, a_v] +
a_u^* [a_s, a_v] a_t \\
& =  
\delta_t(u) a_s a_v +
\delta_s(u) a_t a_v.
\end{align*}
Note that
$$
\lambda^* *  \kappa (t,v)
=
\int_T \lambda^*(t,u) \kappa(u,v) du
=
\int_T \lambda(u,t) \kappa(u,v) du.
$$
Then
\begin{align*}
 [\Xi_{0,2}(\lambda), \Xi_{1,1}(\kappa)]
& =
\int_{T^4} \lambda(s,t) \kappa(u,v) \
[
a_s a_t, a_u^* a_v
]
\ ds dt du dv\\
& = 
\int_{T^3}
\lambda(s,u) \kappa(u,v) a_s a_v ds du  dv+
\int_{T^3}
\lambda(u,t) \kappa(u,v) a_t a_v dt du dv 
\\
& = 
\int_{T^2}
(\lambda *   \kappa)(s,v) a_s a_v ds dv 
+
\int_{T^2}
(\lambda^*  *   \kappa)(t,v) a_t a_v dt dv \\
& =
\Xi_{0,2}(\lambda *   \kappa) 
+ \Xi_{0,2}(\lambda^* *   \kappa). \ \Box
\end{align*}
\begin{lem}
\label{Commutations4}
Let  $\kappa \in E_\mathbb{C} \otimes  E_\mathbb{C}^*$. 
Then
$[\Delta_G, \Xi_{1,1}(\kappa)] = 2\Xi_{0,2}(\kappa)$.
\end{lem}
\noindent
\emph{Proof:}
Recall that $\Delta_G = \Xi_{0,2}(\tau)$. 
Then $\tau^* = \tau$ and
 $\tau * \kappa = \kappa$ since
$\tau$ corresponds to the identity operator in the canonical isomorphism 
$E_\mathbb{C} \otimes E_\mathbb{C}^* \cong
\mathcal{L}(E_\mathbb{C}, E_\mathbb{C})$.
The conclusion follows by Lemma \ref{Commutations3}.
$\Box$ 

\begin{lem}
\label{Commutations5}
Let $\kappa \in E_\mathbb{C} \otimes E_\mathbb{C}^*$ be skew-symmetric.
Then
$
[\Xi_{0,2}(\kappa), \Xi_{1,1}(\kappa)] = 0.
$
\end{lem}
\noindent
\emph{Proof:}
Since $\kappa$ is skew-symmetric,
$\kappa^*(s,t)  = -\kappa(s,t)$
for all $s,t \in T$.
Then
$$
\kappa^* *  \kappa (t,v)
=
\int_T \kappa^*(t,u) \kappa(u,v) du
=
-\int_T \kappa(t,u) \kappa(u,v) du
=
-\kappa * \kappa (t,v).
$$
Applying Lemma \ref{Commutations3}, we have
$$
[\Xi_{0,2}(\kappa), \Xi_{1,1}(\kappa)]
=
\Xi_{0,2}(\kappa * \kappa) 
- \Xi_{0,2} (\kappa * \kappa)
= 0. \ \Box
$$

We record the following commutation relations which are immediate from 
Proposition \ref{generalizedCCR},
Remark \ref{annihilation-commute},
and Lemmas \ref{ManyRelations}, 
\ref{CCR-Number-operator},
\ref{Qwn-derivatives-Xi},
\ref{Commutations4} and
\ref{Commutations5}.

\begin{thm} 
\label{Commutations1}
Let $\zeta \in E_\mathbb{C}$
and $S \in \mathcal{L}(E_\mathbb{C},E_\mathbb{C})$ be skew-symmetric.
Let $S^0 = Id$.
Then for $j, k \in \mathbb{N} \cup \{ 0 \}$, it holds that
\begin{enumerate}
	\item $[a(S^j \zeta), a^*(S^k \zeta)] 
	= \langle S^j \zeta, S^k \zeta \rangle Id$,
	\item $[\Delta_G, \Delta_G(S)] = 0$,
	\item $[N, \Delta_G(S)] = -2 \Delta_G(S)$,
	\item $[N, \Lambda(S)] = 0$,
	\item $[a^*(S^k \zeta) , \Lambda(S)] = a^*(S^{k+1} \zeta)$,
	\item $[a(S^k \zeta) , \Lambda(S)] = (-1)^{k+1} 
	a(S^{k+1} \zeta)$,
	\item $[a^*(S^k\zeta), \Delta_G(S)] = 0$,	 
	\item $[a(S^k\zeta), \Delta_G(S)] = 0$,
	\item $[\Delta_G, \Lambda(S)] = 2\Delta_G(S)$,
	\item $[\Delta_G(S), \Lambda(S)] = 0$.	
\end{enumerate}
\end{thm}




\begin{thm} 
\label{Lie-algebra-result}
For each nonzero $\zeta \in E_\mathbb{C}$
and a skew-symmetric 
$S \in \mathcal{L}( E_\mathbb{C}, E_\mathbb{C})$,
let $S^0 = Id$ and
$$\
\mathfrak{h} = 
\langle Id, a(S^k\zeta), a^*(S^k \zeta), N, \Lambda(S),
 \Delta_G, \Delta_G(S) : k = 0, 1, 2, \dots \rangle 
$$
Then $\mathfrak{h}$ is a (possibly infinite dimensional) non-nilpotent solvable complex Lie algebra.
\end{thm}
\noindent
\emph{Proof:}
From the commutation relations given in
Theorems \ref{ManyRelations} and \ref{Commutations1}
, we see that 
$\mathfrak{h}$ is closed under the Lie bracket.
Moreover, we have
\begin{align*}
\mathfrak{h}^{(1)}
& = [\mathfrak{h}, \mathfrak{h}]
=
\langle
Id, a(S^k \zeta), a^*(S^k \zeta), \Delta_G, \Delta_G(S)
: k = 0, 1, 2, \dots
\rangle
\\
\mathfrak{h}^{(2)}
& = [\mathfrak{h}^{(1)}, \mathfrak{h}^{(1)}]
=
\langle Id, a(S^k \zeta) 
: k = 0,1,2, \dots
\rangle \\
\mathfrak{h}^{(3)}
& = [\mathfrak{h}^{(2)}, \mathfrak{h}^{(2)}]
=
\{0\}.
\end{align*}
Thus $\mathfrak{h}$ is solvable.
Moreover,
$
\mathfrak{h}_{(1)} 
= [\mathfrak{h}, \mathfrak{h}]
$
and
$$
[\mathfrak{h}, \mathfrak{h}_{(1)}]
 = 
\langle Id,  a(S^k\zeta), a^*(S^k\zeta), \Delta_G, \Delta_G(S)
: k = 0, 1, 2, \dots
\rangle
= \mathfrak{h}_{(1)}.
$$
Thus, $\mathfrak{h}$ is non-nilpotent.
$\Box$

\begin{rmk}
It is worth noting that  $S^k$ is skew-symmetric for an odd $k$ while it is symmetric for an even $k$.
The dimension of $\mathfrak{h}$ depends on the properties
of   
$S \in \mathcal{L}(E_\mathbb{C}, E_\mathbb{C})$.
Note that for $\zeta \in E_\mathbb{C}$, 
$\{S^k \zeta\}$ is contained in the the range of $S$.
Therefore, if $S$ is a finite rank operator, i.e.,
$dim \mathcal{R}(S) < \infty$, then
$\mathfrak{h}$ is finite dimensional.
We may also investigate the eigenvectors of $S$.
For instance, adding the condition that $S \zeta = \zeta$
(e.g. Theorem \ref{Fixed-point}) 
makes $\mathfrak{h}$ a seven-dimensional Lie algebra.
More generally, if $\zeta$ is an eigenvector of
$S^k$,
i.e. $S^k \zeta = \lambda \zeta$, for some $\lambda \in \mathbb{C}$,
then $\mathfrak{h}$ is finite dimensional.
By an orbit of $S \in \mathcal{L}(E_\mathbb{C}, E_\mathbb{C})$ we mean a sequence
$\{S^n \xi : n = 0, 1, 2, \dots \}$, where
$\xi \in E_\mathbb{C}$ is a fixed vector.
To read more about orbits of operators,  see 
\cite{Muller2004, Muller2009}.
A vector $e$ is supercyclic for an operator $T$ if the vectors $\xi T^n e$, $\xi$ complex, are dense.
In \cite{Hilden1974}, Hilden and Wallen proved that no normal operator on a Hilbert space of
dimension $> 1$ can be supercyclic.
\end{rmk}

Observe that the one-dimensional subalgebra 
$\langle \Delta_G(S) \rangle$
is an ideal of $\mathfrak{h}$.
We have the following proposition.

\begin{prop} \label{IdealContainsId}
Let $\zeta \in E_\mathbb{C}$ and
$S \in \mathcal{L}(E_\mathbb{C}, E_\mathbb{C})$
be skew-symmetric. 
Set $S^0 = Id$ and
$$
\mathfrak{h} = 
\langle Id, a(S^k\zeta), a^*(S^k \zeta), N, \Lambda(S),
 \Delta_G, \Delta_G(S) : k = 0, 1, 2, \dots \rangle. 
$$
Then every nonzero ideal $\mathfrak{a}$ of $\mathfrak{h}$ 
with $\Delta_G(S) \notin \mathfrak{a}$
contains the identity operator.
\end{prop}
\noindent
\emph{Proof:}
Suppose there exists a nonzero ideal $\mathfrak{a}$ of 
$\mathfrak{h}$ 
with $\Delta_G(S) \notin \mathfrak{a}$
such that $Id \notin\mathfrak{a}$.
Then there exists $x \in \mathfrak{a}$, $x \neq 0$.
\\
\textit{Case 1:} $x = a(S^k \zeta)$. 
Then 
$[a^*(S^k \zeta), a(S^k \zeta)] 
= -\langle S^k \zeta, S^k \zeta \rangle Id$, a contradiction.\\
\textit{Case 2:} $x = a^*(S^k \zeta)$.
Then $[a( S^k\zeta), a^*(S^k\zeta)] 
= \langle S^k \zeta, S^k \zeta \rangle  Id$, a contradiction.\\
\textit{Case 3:} $x = N$.
Note that
$[a(\zeta), N] = a(\zeta) \in \mathfrak{a}$. Then apply Case 1.
\\
\textit{Case 4:} $x = \Lambda(S)$.
Note that
$[a(\zeta), \Lambda(S)] = -a(S\zeta) \in \mathfrak{a}$. Then apply Case 1.
\\
\textit{Case 5:} $x = \Delta_G$.
Note that
$[a^*(\zeta), \Delta_G] = -2a(\zeta) \in \mathfrak{a}$.
Then apply Case 1.
 $\Box$

\begin{thm} \label{NotSemisimple}
The complex Lie algebra 
$\mathfrak{h}
$ 
described in Theorem \ref{Lie-algebra-result}
is not semisimple. 
\end{thm}
\noindent
\emph{Proof:}
The proof of Theorem \ref{Lie-algebra-result}
shows that $\mathfrak{h}$ contains a nonzero solvable ideal, e.g. $\mathfrak{h}^{(1)}$. 
$\Box$


\end{document}